\documentclass[a4paper,11pt] {amsart}
\usepackage{amsthm,amssymb,latexsym,mathrsfs}
\input amssym.def

\author{Carnot D. Kenfack and Beno\^it F. Sehba}

\title[Maximal function and Carleson measures ]{Maximal function and Carleson measures in B\'ekoll\'e-Bonami weights}
\newtheorem{theorem}{T{\hskip 0pt\footnotesize\bf HEOREM}}[section]
\newtheorem{lemma}[theorem]{L{\hskip 0pt\footnotesize\bf EMMA}}
\newtheorem{proposition}[theorem]{P{\hskip 0pt\footnotesize\bf ROPOSITION}}

\newtheorem{corollary}[theorem]{C{\hskip 0pt\footnotesize\bf OROLLARY}}

\newcommand{\bprop} {\begin{proposition}}
\newcommand{\eprop} {\end{proposition}}
\newcommand{\btheo} {\begin{theorem}}
\newcommand{\etheo} {\end{theorem}}
\newcommand{\blem} {\begin{lemma}}
\newcommand{\elem} {\end{lemma}}
\newcommand{\bcor} {\begin{corollary}}
\newcommand{\ecor} {\end{corollary}}

\newcommand{\Be}{\begin{equation}}
\newcommand{\Ee}{\end{equation}}
\newcommand{\Bea}{\begin{eqnarray}}
\newcommand{\Eea}{\end{eqnarray}}
\newcommand{\Bes}{\begin{equation*}}
\newcommand{\Ees}{\end{equation*}}
\newcommand{\Beas}{\begin{eqnarray*}}
\newcommand{\Eeas}{\end{eqnarray*}}
\newcommand{\Ba}{\begin{array}}
\newcommand{\Ea}{\end{array}}

\scrollmode

\begin{document}
\address{Carnot D. Kenfack, D\'epartement de math\'ematiques, Facult\'e des Sciences, Universit\'e de Yaound\'e I, BP.812 Yaound\'e, Cameroun}
\email{dondjiocarnot@yahoo.fr}
\address{Beno\^it F. Sehba, }
\email{bsehba@gmail.com}


\begin{abstract}
Let $\omega$ be a B\'ekoll\'e-Bonami weight. We give a complete characterization of the positive measures $\mu$ such that
$$\int_{\mathcal H}|M_\omega f(z)|^qd\mu(z)\le C\left(\int_{\mathcal H}|f(z)|^p\omega(z)dV(z)\right)^{q/p}$$
and
$$\mu\left(\{z\in \mathcal H: Mf(z)>\lambda\}\right)\le \frac{C}{\lambda^q}\left(\int_{\mathcal H}|f(z)|^p\omega(z)dV(z)\right)^{q/p}$$
where $M_\omega$ is the weighted Hardy-Littlewood maximal function on the upper-half plane $\mathcal H$, and $1\le p,q<\infty$.
\end{abstract}
\maketitle
\section{Introduction}
Let $\mathcal H$ be the upper-half plane, that is the set $\{z=x+iy\in \mathbb C: x\in \mathbb R,\,\,\, \textrm{and}\,\,\,y>0\}$. Given $\omega$ a nonnegative locally integrable function on $\mathcal H$ (i.e a weight), and $1\le p<\infty$, we denote by $L_\omega^p(\mathcal H)$, the set of functions $f$ defined on $\mathcal H$ such that
$$||f||_{p,\omega }^p:=\int_{\mathbb D}|f(z)|^p\omega(z)dV(z)<\infty$$
with $dV$ being the Lebesgue measure on $\mathcal H$.
\vskip .1cm
Given a weight $\omega$, and $1<p<\infty$, we say $\omega$ is in the B\'ekoll\'e-Bonami class $B_p$,  if
$$[\omega]_{B_p}:=\sup_{I\subset \mathbb R,\,\,\, I\,\,\,\textrm{interval} }\left(\frac{1}{|Q_I|}\int_{Q_I}\omega(z)dV(z)\right)\left(\frac{1}{|Q_I|}\int_{Q_I}\omega(z)^{1-p'}dV(z)\right)^{p-1}<\infty,$$
$Q_I:=\{z=x+iy\in \mathbb C:x\in I\,\,\,\textrm{and}\,\,\,0<y<|I|\}$, $|Q_I|=\int_{Q_I}dV(z)$, $pp'=p+p'$. This is the exact range of weights $\omega$ for which the orthogonal projection $P$ from $L^2(\mathcal H, dV(z))$ to its closed subspace consisting of analytic functions is bounded on $L_\omega^p(\mathcal H)$ (see \cite{Bek,BB, PR, Sehba}).
\vskip .1cm
Let $1<p<\infty$, and $\omega\in B_p$. We provide in this note a full characterization of positive measures $\mu$ on $\mathcal H$ such that the following Carleson-type embedding
\Be\label{eq:question1}
\int_{\mathcal H}|M_\omega f(z)|^qd\mu(z)\le C\left(\int_{\mathcal H}|f(z)|^p\omega(z)dV(z)\right)^{q/p}
\Ee
holds when $p\le q<\infty$ and when $p>q$, where $M_\omega$ is the weighted Hardy-Littlewood maximal function, $$M_\omega f(z):=\sup_{ I\,\,\,\textrm{interval in} \,\,\,\mathbb R,\,\, z\in Q_I}\frac{1}{|Q_I|_\omega}\int_{Q_I}|f(z)|\omega(z)dV(z),$$ $|Q_I|_\omega=\omega(Q_I)=\int_{Q_I}\omega(z)dV(z)$.
\vskip .1cm
We also characterize those positive measures $\mu$ on $\mathcal H$ such that
\Be\label{eq:question2}
\mu\left(\{z\in \mathcal H: Mf(z)>\lambda\}\right)\le \frac{C}{\lambda^q}\left(\int_{\mathcal H}|f(z)|^p\omega(z)dV(z)\right)^{q/p}
\Ee
where $M$ is the unweighted Hardy-Littlewood maximal function ($M=M_\omega$ with $\omega(z)=1$ for all $z\in \mathcal H$).
\vskip .1cm
Before stating our main results, let us see how the above questions are related to some others in complex analysis. We recall that the Bergman space $A_\omega^p(\mathcal H)$ is the subspace of $L_\omega^p(\mathcal H)$ consisting of holomorphic functions on $\mathcal H$. The usual Bergman spaces in the unit disc of $\mathbb C$ or the unit ball of $\mathbb C^n$ correspond to the weights $\omega(z)=(1-|z|^2)^\alpha dA(z)$, $\alpha>-1$. A positive measure on $\mathcal H$ is called a $q$-Carleson measure for $A_\omega^p(\mathcal H)$ if there is a constant $C>0$ such that for any $f\in A_\omega^p(\mathcal H)$,
\Be\label{eq:Carlmeadef}
\int_{\mathcal H}|f(z)|^qd\mu(z)\le C\left(\int_{\mathcal H}|f(z)|^p\omega(z)dV(z)\right)^{q/p}.
\Ee
Carleson measures are very useful in the study of many other questions in complex and harmonic analysis: Toeplitz operators, Ces\`aro-type integrals, embeddings between different analytic function spaces, etc... Carleson measures for Bergman spaces with standard weights $\omega(z)=(1-|z|^2)^\alpha dA(z)$ in the unit disc and the unit ball of $\mathbb C^n$, $\alpha>-1$ have been studied in \cite{CW, Hastings, Luecking1, Luecking3, Oleinik, stegenga}. The case of Bergman spaces of the unit disc of $\mathbb C$ with B\'ekoll\'e-Bonami weights has been handled in \cite{Constantin, Gu}.
\vskip .1cm
Let us suppose that $\omega\in B_p$. Applying the mean value property one obtains that there is a constant $C>0$ such that for any $f\in A_\omega^p(\mathcal H)$,  and for any $z\in \mathcal H$, if $I$ is the unique interval such that $Q_I$ is centered at $z$, then
$$|f(z)|\le \frac{C}{\omega (Q_I)}\int_{Q_I}|f(w)|\omega(w)dV(w).$$
It follows that any measure satisfying (\ref{eq:question1}) is a $q$-Carleson measure for $A_\omega^p(\mathcal H)$.
\vskip .1cm
Our first main result is the following.
\btheo\label{thm:main1}
Let $1<p\le q<\infty$, and $\omega$ a weight on $\mathcal H$. Assume that $\omega\,\in B_{p}$. Then the following assertions are equivalent.
 \begin{itemize}
\item[(1)] There exists a constant $C_1>0$ such that for any $f\in L_\omega^p(\mathcal H)$,
\Be\label{eq:main11}
\left(\int_{\mathcal H}|M_\omega f(z)|^qd\mu(z)\right)^{1/q}\le C_1\|f\|_{p,\omega}.
\Ee

\item[(2)] There is a constant $C_2$ such that for any interval $I\subset \mathbb R$,
\Be\label{eq:main12} \mu(Q_I)\leq C_2 (\omega(Q_I))^{\frac{q}{p}}.\Ee
\end{itemize}
\etheo
Our next result provides estimations with loss.
\btheo\label{thm:main2}
Let $1<q<p<\infty$, and $\omega$ a weight on $\mathcal H$. Assume that $\omega\,\in B_{p}$. Then (\ref{eq:main11}) holds if and only if the function
\Be\label{eq:main21} K_{\mu}(z):=\sup_{I\subset \mathbb R,\,\,\,I\,\,\,\textrm{interval},\,\,\,z\in Q_I}\frac{\mu(Q_{I})}{\omega(Q_{I})}
\Ee belongs to
$L_{\omega}^s(\mathcal H)$ where $s=\frac{p}{p-q}$.
\etheo
Our last result provides weak-type estimates.
\btheo\label{thm:main3}
Let $1\le p, q<\infty$, and $\omega$ a weight on $\mathcal H$. Then the following assertions are equivalent.
\begin{itemize}
\item[(a)] There is a constant $C_1>0$ such that for any $f\in L_\omega^p(\mathcal H)$, and any $\lambda>0$,
\Be\label{eq:main31}
\mu\left(\{z\in \mathcal H: Mf(z)>\lambda\}\right)\le \frac{C_1}{\lambda^q}\left(\int_{\mathcal H}|f(z)|^p\omega(z)dV(z)\right)^{q/p}
\Ee
\item[(b)] There is a constant $C_2>0$ such that for any interval $I\subset \mathbb R$,
\Be\label{eq:main32}
|Q_I|^{-q/p}\left(\frac{1}{|Q_I|}\int_{Q_I}\omega^{1-p'}(z)dV(z)\right)^{q/p'}\mu(Q_I)\le C_1
\Ee
where $\left(\frac{1}{|Q_I|}\int_{Q_I}\omega^{1-p'}(z)dV(z)\right)^{1/p'}$ is understood as $\left(\inf_{Q_I}\omega\right)^{-1}$ when $p=1$.
\item[(c)] There exists a constant $C_3>0$ such that for any locally integrable function $f$ and  any interval $I\subset \mathbb R$,
\Be\label{eq:main33}
\left(\frac{1}{|Q_I|}\int_{Q_I}|f(z)|dV(z)\right)^q\mu(Q_I)\le C_3\left(\int_{Q_I}|f(z)|^p\omega(z) dV(z)\right)^{q/p}.
\Ee
\end{itemize}
\etheo
A special case of Theorem \ref{thm:main3} appears when $\mu$ is a continuous measure with respect to the Lebesgue measure $dV$ in this sense that $d\mu(z)=\sigma(z)dV(z)$, this provides a weak-type two-weight norm inequality for the maximal function.
\vskip .1cm
To prove the sufficient part in the three theorems above, we will observe that the matter can be reduced to the case of the dyadic maximal function. We then use an idea that comes from real harmonic analysis (see for example \cite{Cruz1,cuervafrancia,Sawyer}) and consists of discretizing  integrals using appropriate level sets and in our case, the nice properties of the upper-halfs of Carleson boxes when they are supported by dyadic intervals. For the proof of the necessity in Theorem \ref{thm:main2}, let us observe that when comes to estimations with loss for the case of the usual Carleson measures for analytic functions, one needs atomic decomposition of functions in the Bergman spaces to apply a method developed by D. Luecking \cite{Luecking1}. We do not see how this can be extended here and instead, we show that one can restrict to the dyadic case, and use boundedness of the maximal functions and a duality argument. We note that a duality argument has been used for the same-type of question for weighted Hardy spaces in \cite{Gu1}.
\vskip .1cm
 Given two positive quantities $A$ and $B$, the notation $A\lesssim B$ (resp. $B\lesssim A$) will mean that there is an universal constant $C>0$ such that $A\le CB$ (resp. $B\le CA$).
\section{Useful observations and results}
Given an interval $I\subset \mathbb R$, the upper-half of the Carleson box $Q_I$ associated to $I$ is the subset $T_I$ defined by $$T_I:=\{z=x+iy\in \mathbb C: x\in I,\,\,\,\textrm{and}\,\,\,\frac{|I|}{2}<y<|I|\}.$$
Note that $|Q_I|\backsimeq |T_I|$. We observe the following weighted inequality.
\blem\label{lem:ineqomegaQT}
Let $1<p<\infty$. Assume that $\omega$ belongs to the
B\'{e}koll\'{e}-Bonami class $B_p$. Then there is a constant $C>0$ such that for any interval $I\subset \mathbb R$,
$$\omega(Q_I)\le C
[\omega]_{B_p}\omega(T_I).$$
\elem
\begin{proof}
Using H\"older's inequality and the definition of B\'ekoll\'e-Bonami weight, we obtain
\Beas
\frac{|T_I|^p}{|Q_I|^p} &\le& \frac{1}{|Q_I|^p}\left(\int_{T_I}\omega(z)dV(z)\right)\left(\int_{T_I}\omega^{-p'/p}(z)dV(z)\right)^{p/p'}\\ &\le& \frac{1}{|Q_I|^p}\left(\int_{T_I}\omega(z)dV(z)\right)\left(\int_{Q_I}\omega^{-p'/p}(z)dV(z)\right)^{p/p'}\\ &\le& [\omega]_{B_p}\frac{\omega(T_I)}{\omega(Q_I)}.
\Eeas
Thus $\omega(Q_I)\le [\omega]_{B_p}\left(\frac{|Q_I|}{|T_I|}\right)^p\omega(T_I)\backsimeq [\omega]_{B_p}\omega(T_I).$
\end{proof}
We will also need the following lemma.
\blem\label{lem:equivdefBpq}
Let  $1\le p, q <\infty$ and suppose that $\omega$ is a weight, and $\mu$ a positive measure on $\mathcal H$. Then the following assertions are equivalent.
\begin{itemize}
\item[(i)] There exists a constant $C_1>0$ such that for any interval $I\subset \mathbb R$,
\Be\label{eq:equivMuBpq1}
|Q_I|^{-q/p}\left(\frac{1}{|Q_I|}\int_{Q_I}\omega^{1-p'}(z)dV(z)\right)^{q/p'}\mu(Q_I)\le C_1
\Ee
where $\left(\frac{1}{|Q_I|}\int_{Q_I}\omega^{1-p'}(z)dV(z)\right)^{1/p'}$ is understood as $\left(\inf_{Q_I}\omega\right)^{-1}$ when $p=1$.
\item[(ii)] There exists a constant $C_2>0$ such that for any locally integrable function $f$ and  any interval $I\subset \mathbb R$,
\Be\label{eq:equivMuBpq2}
\left(\frac{1}{|Q_I|}\int_{Q_I}|f(z)|dV(z)\right)^q\mu(Q_I)\le C_2\left(\int_{Q_I}|f(z)|^p\omega(z) dV(z)\right)^{q/p}.
\Ee
\end{itemize}
\elem
\begin{proof}
That $(\textrm{ii})\Rightarrow (\textrm{i})$ follows by testing $(\textrm{ii})$ with $f(z)=\chi_{Q_I}(z)\omega^{1-p'}(z)$ if $p>1$. For $p=1$, take $f(z)=\chi_{S}(z)$ where $S$ a subset of $Q_I$. One obtains that $$\frac{\mu(Q_I)}{|Q_I|^q}\le C_2\left(\frac{\omega(S)}{|S|}\right)^q.$$ As this happens for any subset $S$ of $Q_I$, it follows that for any $z\in Q_I$, $$\frac{\mu(Q_I)}{|Q_I|^q}\le C_2\left(\omega(z)\right)^q$$ which implies (\ref{eq:equivMuBpq1}) for $p=1$.
 \vskip .1cm
 Let us check that $(\textrm{i})\Rightarrow (\textrm{ii})$. Applying H\"older's inequality (in case $p>1$) to the right hand side of (\ref{eq:equivMuBpq2}), we obtain
\Beas
&& \left(\frac{1}{|Q_I|}\int_{Q_I}|f(z)|dV(z)\right)^q\mu(Q_I)\\ &\le& |Q_I|^{-q}\left(\int_{Q_I}\omega^{-p'/p}(z) dV(z)\right)^{q/p'}\mu(Q_I)\left(\int_{Q_I}|f(z)|^p\omega(z) dV(z)\right)^{q/p}\\ &\le& C\left(\int_{Q_I}|f(z)|^p\omega(z) dV(z)\right)^{q/p}.
\Eeas
For $p=1$, we easily obtain
\Beas
\left(\frac{1}{|Q_I|}\int_{Q_I}|f(z)|dV(z)\right)^q\mu(Q_I) &\le& \frac{\left(\inf_{Q_I}\omega\right)^{-q}}{|Q_I|^q}\left(\int_{Q_I}|f(z)|\omega(z)dV(z)\right)^q\mu(Q_I)\\ &\le& C\left(\int_{Q_I}|f(z)|\omega(z)dV(z)\right)^q.
\Eeas
The proof is complete.
\end{proof}

Next, we consider the following system of dyadic grids,
$$\mathcal D^\beta:=\{2^j\left([0,1)+m+(-1)^j\beta\right):m\in \mathbb Z,\,\,\,j\in \mathbb Z \},\,\,\,\textrm{for}\,\,\,\beta\in \{0,1/3\}.$$
For more on this system of dyadic grids and its applications, we refer to \cite{AlPottReg, HyLaPerez, HyPerez, Lerner,LerOmbroPerezetal, PR, Sehba}. When $\beta=0$, we use the notation $\mathcal D=\mathcal D^0$ that we call the standard dyadic grid of $\mathbb R$. When $I$ is a dyadic interval, we denote by $I^-$ and $I^+$ its left half and its right half respectively. We make the following observation which is surely known.
\blem\label{lem:dyacovering}
Any interval $I$ of $\mathbb R$ can be covered by at most two adjacent dyadic intervals $I_1$ and $I_2$ in the same dyadic grid such that
$$|I|<|I_1|=|I_2|\le 2|I|.$$
\elem
\begin{proof}
Without loss of generality, we can suppose that $I=[a,b)$. For $x\in \mathbb R$, we denote by $[x]$ the unique integer such that $[x]\le x<[x]+1$. If $I\in \mathcal D$, then there is nothing to say. If $|I|=1$, then the dyadic interval $[k,k+1)$ where $k=[a]$ covers $I$.
\vskip .1cm
Let us suppose in general that $I$ is not dyadic. Let $j$ be the unique integer such that \Be\label{eq:embedinter} 2^{-j}\le b-a=|I|<2^{-j+1},\Ee and define the set $$E_{a,b}:=\{l\in \mathbb Z:a<l2^{-j}\le b\}.$$ Then $E_{a,b}$ is not empty. To see this, take $k=[a2^j]$, then $(k+1)2^{-j}\le b$ since if not, we will have $[a,b)\subset [k2^{-j},(k+1)2^{-j})$ and consequently, $|I|=b-a<|[k2^{-j},(k+1)2^{-j})|=2^{-j}$ which contradicts (\ref{eq:embedinter}). Let $$k_0:=\max\{k:k\in E_{a,b}\}.$$ Then we necessarily have $(k_0-2)2^{-j}\le a$ since if not, $|I|=b-a>k_02^{-j}-a>2^{-j+1}$ and this contradicts (\ref{eq:embedinter}).

As from the definition of $k_0$ we have $b\le (k_0+1)2^{-j}$, it comes that if $(k_0-1)2^{-j}\le a$, then the union $[(k_0-1)2^{-j},k_02^{-j})\cup [k_02^{-j},(k_0+1)2^{-j})$ covers I, and taking $I_1$ and $I_2$ such that $I_1^+=[(k_0-1)2^{-j},k_02^{-j})$ and $I_2^-=[k_02^{-j},(k_0+1)2^{-j})$ we get the lemma. If $(k_0-1)2^{-j}>a$, then $I\subset I_1\cup I_2$ where $I_1=[(l_0-1)2^{-j+1},l_02^{-j+1})$, $I_2=[l_02^{-j+1},(l_0+1)2^{-j+1})$ with $k_0=2l_0$ if $k_0$ is even or $k_0=2l_0+1$ otherwise.
The proof is complete.

\end{proof}
\section{Proof of the results}
Let us start with some observations.  Recall that given $Q_I$, its upper-half is  the set $$T_I:=\{x+iy\in \mathcal H: x\in I,\,\,\,\textrm{and}\,\,\,\frac{|I|}{2}<y<|I|\}.$$ It is clear that the family $\{T_I\}_{I\in \mathcal D}$ where $\mathcal D$ is a dyadic grid in $\mathbb R$ provides a tiling of $\mathcal H$.
\vskip .1cm
Next we recall with \cite{PR} that given an interval $I\subset \mathbb R$, there is a dyadic interval $K\in \mathcal D^\beta$ for some $\beta\in \{0,1/3\}$ such that $I\subseteq K$ and $|K|\le 6|I|$. It follows in particular that $|Q_K|\le 36|Q_I|$. Also, proceeding as in the proof of Lemma \ref{lem:ineqomegaQT} one obtains that $\omega(Q_K)\lesssim [\omega]_{B_p}\omega(Q_I)$. It follows that
$$\frac{1}{\omega(Q_I)}\int_{Q_I}|f(z)|\omega(z)dV(z)\lesssim \frac{1}{\omega(Q_K)}\int_{Q_K}|f(z)|\omega(z)dV(z)$$
and consequently that for any locally integrable function $f$,
\Be\label{eq:Maxfunctdyaineq}
M_\omega f(z)\lesssim \sum_{\beta\in \{0,1/3\}}M_{d,\omega}^\beta f(z),\,\,\,z\in \mathcal H
\Ee
where $M_{d,\omega}^\beta$ is defined as $M_\omega$ but with the supremum taken only over dyadic intervals of the dyadic grid $\mathcal D^\beta$. When $\omega\equiv 1$, we use the notation $M_d^\beta$, and if moreover, $\beta=0$, we just write $M_d$. In the sequel, we will be proving anything only for the case $\beta=0$ which is enough and in this case, we write everything without the superscript $\beta=0$.
\subsection{Proof of Theorem \ref{thm:main1}}
First suppose that (\ref{eq:main11}) holds and observe that for any interval $I\subset \mathbb R$, $1\le M_\omega \chi_{Q_I}(z)$ for any $z\in Q_I$. It follows that
$$\left(\mu(Q_I)\right)^{1/q}\le \left(\int_{\mathcal H}\left(M_\omega \chi_{Q_I}(z)\right)^qd\mu(z)\right)^{1/q}\le C_1\|\chi_{Q_I}\|_{p,\omega}=\left(\omega(Q_I)\right)^{1/p}$$
which provides that for any interval $I\subset \mathbb R$,
$$\mu(Q_I)\le C_1\left(\omega(Q_I)\right)^{q/p}.$$
That is (\ref{eq:main12}) holds.

To prove that $(\textrm{ii})\Rightarrow (\textrm{i})$, it is enough by the observations made at the beginning of this section to prove the following.
\blem\label{lem:main1}
Let $1< p\le q<\infty$. Assume that $\omega$ is a weight in the class $B_p$ such that (\ref{eq:main12}) holds. Then there is a positive constant $C$ such that for any $f\in L_\omega^p(\mathcal H)$,
\Be\label{eq:main13}
\left(\int_{\mathcal H}|M_{d,\omega} f(z)|^qd\mu(z)\right)^{1/q}\le C_1\|f\|_{p,\omega}.
\Ee
\elem
\begin{proof}
Let $a\ge 2$. To each integer $k$, we associate the set
$$\Omega_{k}:=\{z\in \mathcal H: a^k<M_{d,\omega}f(z)\leq a^{k+1}\}.$$
We observe that
$\Omega_{k}\subset \cup_{j=1}^{\infty}Q_{I_{k,j}},$ where
$Q_{I_{k,j}}$ is a dyadic cube maximal (with respect to the inclusion) such that
$$\frac{1}{\omega(Q_{I_{k,j}})}\int_{Q_{I_{k,j}}}|f(z)|\omega(z)dV(z)>a^k.$$
It follows using Lemma \ref{lem:ineqomegaQT} that
 \Beas
 \int_{\mathcal H}(M_{d,\omega}f(z))^{q}d\mu(z) &=& \sum_{k}\int_{\Omega_k}(M_{d,\omega}f(z))^{q}d\mu(z)\\
 &\le& a^{q}\sum_{k}a^{kq}\mu(\Omega_k)\\
 &\le& a^{q}\sum_{k,j}a^{kq}\mu(Q_{I_{k,j}})\\
 &\lesssim&
 a^{q}\sum_{k,j}\left(\frac{1}{\omega(Q_{I_{k,j}})}\int_{Q_{I_{k,j}}}|f(z)|\omega(z)dV(z)\right)^{q}\mu(Q_{I_{k,j}})\\
 &\lesssim&
 a^{q}\sum_{k,j}\left(\frac{1}{\omega(Q_{I_{k,j}})}\int_{Q_{I_{k,j}}}|f(z)|\omega(z)dV(z)\right)^{q}(\omega(Q_{I_{k,j}}))^{\frac{q}{p}}\\
&\lesssim&
a^{q}\left(\sum_{k,j}\left(\frac{1}{\omega(Q_{I_{k,j}})}\int_{Q_{I_{k,j}}}|f(z)|\omega(z)
dV(z)\right)^p\omega(Q_{I_{k,j}})\right)^{\frac{q}{p}}\\
 &\lesssim&
 a^{q}\left(\sum_{k,j}\left(\frac{1}{\omega(Q_{I_{k,j}})}\int_{Q_{I_{k,j}}}|f(z)|\omega(z)
 dA(z)\right)^p[\omega]_{B_p}\omega(T_{I_{k,j}})\right)^{\frac{q}{p}}\\
&\lesssim& \left(\sum_{k,j}\int_{T_{I_{k,j}}}\left(\frac{1}{\omega(Q_{I_{k,j}})}\int_{Q_{I_{k,j}}}|f(z)|\omega(z)dV(z)\right)^p\omega(w)dV(w)\right)^{q/p}\\
&\lesssim& \left(\sum_{k,j}\int_{T_{I_{k,j}}}(M_{d,\omega}f(z))^p
dV(z)\right)^{\frac{q}{p}}\\
&\lesssim& \left(\int_{\mathcal H}|f(z)|^p\omega(z)
dV(z)\right)^{\frac{q}{p}}.
\Eeas
The proof of the lemma is complete.
\end{proof}
\subsection{Proof of Theorem \ref{thm:main2}}
Let us start by proving the following lemma.
%
\blem\label{lem:main21}
Let $1\le q<p<\infty$, and let $\omega$ be a weight in the class $B_p$. Assume that $\mu$ is a positive measure on $\mathcal H$ such that the function $K_\mu$ defined by (\ref{eq:main21}) belongs to $L_\omega^s(\mathcal H)$, $s=\frac{p}{p-q}$. Then there is a constant $C_1>0$ such that for any $f\in L_\omega^p(\mathcal H)$, (\ref{eq:main13}) holds.
\elem
\begin{proof}
We proceed as in the proof of Lemma \ref{lem:main1}, using the same notations.
Using H\"older's inequality and Lemma \ref{lem:ineqomegaQT}, we obtain
 \Beas
 \int_{\mathcal H}(M_{d,\omega}f(z))^{q}d\mu(z) &=& \sum_{k}\int_{\Omega_k}(M_{d,\omega}f(z))^{q}d\mu(z)\\
 &\le& a^{q}\sum_{k}a^{kq}\mu(\Omega_k)\\
 &\le& a^{q}\sum_{k,j}a^{kq}\mu(Q_{I_{k,j}})\\
 &\lesssim&
 a^{q}\sum_{k,j}\left(\frac{1}{\omega(Q_{I_{k,j}})}\int_{Q_{I_{k,j}}}|f(z)|\omega(z)dV(z)\right)^{q}\mu(Q_{I_{k,j}})\\
 &=& a^{q}\sum_{k,j}\left(\frac{1}{\omega(Q_{I_{k,j}})}\int_{Q_{I_{k,j}}}|f(z)|\omega(z)dV(z)\right)^{q}
 \frac{\mu(Q_{I_{k,j}})}{\omega(Q_{I_{k,j}})}\omega(Q_{I_{k,j}})\\
&\lesssim& A^{q/p}B^{1/s}
\Eeas
where $$A=\sum_{k,j}\left(\frac{1}{\omega(Q_{I_{k,j}})}\int_{Q_{I_{k,j}}}|f(z)|\omega(z)dA(z)\right)^{p}\omega(Q_{I_{k,j}})$$
and
$$B=\sum_{k,j}\left(\frac{\mu(Q_{I_{k,j}})}{\omega(Q_{I_{k,j}})}\right)^{s}\omega(Q_{I_{k,j}}).$$
 From the proof of Lemma \ref{lem:main1}, we already know how to estimate $A$. Let us estimate $B$.
\Beas
B &:=&\left(\sum_{k,j}\left(\frac{\mu(Q_{I_{k,j}})}{\omega(Q_{I_{k,j}})}\right)^{s}\omega(Q_{I_{k,j}})\right)\\
&\lesssim&
\sum_{k,j}\left(\frac{\mu(Q_{I_{k,j}})}{\omega(Q_{I_{k,j}})}\right)^{s}\omega(T_{I_{k,j}})\\
&\lesssim& \sum_{k,j}\int_{T_{I_{k,j}}}\left(\frac{\mu(Q_{I_{k,j}})}{\omega(Q_{I_{k,j}})}\right)^{s}\omega(z)dV(z)\\
&\lesssim& \sum_{k,j}\int_{T_{I_{k,j}}}\left(K_{\mu}(z)\right)^{s}\omega(z)dV(z)\\
&\lesssim&
\int_{\mathcal H}\left(K_{\mu}(z)\right)^{s}\omega(z)dV(z)=\|K_{\mu}\|_{s,\omega}^s.
\Eeas

The proof of the lemma is complete.
\end{proof}
We can now prove the theorem
\begin{proof}[Proof of Theorem \ref{thm:main2}]
The proof of the sufficiency follows from Lemma \ref{lem:main21} and the observations made at the beginning of this section. Let us prove the necessity. For this we do the following observations: first, that condition (\ref{eq:main11}) implies that there exists a constant $C>0$ such that for any $f\in L_\omega^p(\mathcal H),$
\Be\label{eq:main111}
\int_{\mathcal H}\left(M_{d,\omega}^\beta f(z)\right)^qd\mu(z)\le C\|f\|_{p,\omega}^q.
\Ee
Second, writing $$K_{d,\mu}^\beta(z):=\sup_{I\in \mathcal {D}^\beta, z\in Q_I}\frac{\mu(Q_I)}{\omega(Q_I)},$$
it is easy to see that for any $z\in \mathcal H$, $$K_\mu(z)\lesssim \sum_{\beta\in \{0,\frac{1}{3}\}}K_{d,\mu}^\beta(z).$$
Thus to prove that $K_\mu\in L_\omega^s(\mathcal H)$ if (\ref{eq:main11}) holds, it is enough to prove that (\ref{eq:main111})
implies that $K_{d,\mu}^\beta\in L_\omega^s(\mathcal H)$. We do this for the standard dyadic grid, i.e for $\beta=0$.
\vskip .1cm
 For $z\in \mathcal H$, we write $Q_z=Q_{I_z}$ ($I_z\in \mathcal D$) for the smallest Carleson box containing $z$, and consider the following weighted box kernel $$K_{d,\omega}(z_0,z):=\frac{1}{\omega(Q_{z_0})}\chi_{Q_{z_0}}(z).$$
For $f$ a locally integrable function, we define
$$K_{d,\omega} f(z_0)=\int_{\mathcal H}K_{d,\omega}(z_0,z)f(z)\omega(z)dV(z)=\frac{1}{\omega(Q_{z_0})}\int_{Q_{z_0}}f(z)\omega(z)dV(z).$$
Finally, we define a function $g$ on $\mathcal H$ by
$$g(z):=\int_{\mathcal H}K_{d,\omega}(\xi,z)d\mu(\xi)=\int_{\mathcal H} \frac{\chi_{Q_\xi}(z)}{\omega(Q_\xi)}d\mu(\xi).$$
For any (dyadic) Carleson box $Q_I$, $I\in \mathcal D$,  writing $Q$ for $Q_I$, we obtain
\Beas
\frac{1}{\omega(Q)}\int_Q g(z)\omega(z)dV(z) &=& \frac{1}{\omega(Q)}\int_Q \left(\int_{\mathcal H} K_{d,\omega}(w,z)d\mu(w)\right)\omega(z)dV(z)\\ &=& \int_{\mathcal H}\int_{\mathcal H}\frac{1}{\omega(Q)}\frac{\chi_{Q_w}(z)\chi_Q(z)}{\omega(Q_w)}\omega(z)dV(z)d\mu(w)\\ &\ge& \int_{Q}\frac{1}{\omega(Q)}\int_{\mathcal H}\frac{\chi_{Q_w\cap Q}(z)}{\omega(Q_w)}\omega(z)dV(z)d\mu(w)\\ &\gtrsim& \frac{1}{\omega(Q)}\int_{Q}d\mu(w)=\frac{\mu(Q)}{\omega(Q)}
\Eeas
Thus for any $z\in \mathcal H$,
$$M_{d,\omega} g(z)\gtrsim \sup_{I\in \mathcal D,\,\,\,\,\,z\in Q_I}\frac{\mu(Q_I)}{\omega(Q_I)}:=K_{d,\mu}(z).$$
Hence if the function $g$ belongs to $L_\omega^s(\mathcal H)$, then
$$\|K_{d,\mu}\|_{s,\omega}\lesssim \|M_{d,\omega} g\|_{s,\omega}\lesssim \|g\|_{s,\omega}.$$
To finish the proof, we only need to check that $g\in L_\omega^s(\mathcal H)$ whenever (\ref{eq:main111}) holds.
 \vskip .1cm
 Let us start by observing the following inequality  between $K_{d,\omega} f$ and $M_{d,\omega} f$. Let $z_0$ be fixed in $\mathcal H$. For any $\xi\in Q_{z_0}$, we have $$K_{d,\omega} f(z_0):=\frac{1}{\omega(Q_{z_0})}\int_{Q_{z_0}}f(z) \omega(z)dV(z)\le M_{d,\omega} f(\xi).$$
 Thus
 \Be\label{eq:KMineq}
 |K_{d,\omega} f(z)|^{1/q}\le M_{d,\omega}\left((M_{d,\omega} f)^{1/q}\right)(z),\,\,\,\textrm{for any}\,\,\,z\in \mathcal H.
 \Ee
 Now, for any $f\in L_\omega^{p/q}(\mathcal H)$, using (\ref{eq:KMineq}), (\ref{eq:main111}) and the boundedness of the maximal function, we obtain
 \Beas
 \left|\int_{\mathcal H}g(z)f(z)\omega(z)dV(z)\right| &=& \left|\int_{\mathcal H}\left(\int_{\mathcal H}K_{d,\omega}(\xi,z)d\mu(\xi)\right)f(z)\omega(z)dV(z)\right|\\ &=& \left|\int_{\mathcal H}\left(\int_{\mathcal H}K_{d,\omega}(\xi,z)f(z)\omega(z)dV(z)\right)d\mu(\xi)\right|\\ &=&
 \left|\int_{\mathcal H}K_{d,\mu}f(\xi)d\mu(\xi)\right|\\ &\le& \int_{\mathcal H}|K_{d,\mu}f(\xi)|d\mu(\xi)\\ &=&
 \int_{\mathcal H}(|K_{d,\mu}f(\xi)|^{1/q})^qd\mu(\xi)\\ &\lesssim& \int_{\mathcal H}\left(M_{d,\omega}\left((M_{d,\omega} f)^{1/q}\right)(\xi)\right)^qd\mu(\xi)\\ &\lesssim&
 \left(\int_{\mathcal H}\left(M_{d,\omega} f(z)\right)^{p/q}\omega(z)dV(z)\right)^{q/p}\\
 &\lesssim& \left(\int_{\mathcal H}|f(z)|^{p/q}\omega(z)dV(z)\right)^{q/p}.
 \Eeas
 Thus there is a constant $C>0$ such that
 $$\|g\|_{s,\omega}:=\sup_{f\in L_\omega^{p/q}(\mathcal H), \|f\|_{p/q,\omega}\le 1}\left|\int_{\mathcal H}g(z)f(z)\omega(z)dV(z)\right|\le C.$$
 The proof is complete.
\end{proof}
\subsection{Proof of Theorem \ref{thm:main3}}
 We start by the following lemma which tells us that we will only need to restrict to level sets involving the dyadic maximal function.
\blem\label{lem:levelsetsembed}
Let $f$ be a locally integrable function. Then for any $\lambda>0$,
\Be\label{eq:levelsetsembed}
\{z\in \mathcal H: Mf(z)>\lambda\}\subset \{z\in \mathbb D: M_df(z)>\frac{\lambda}{68}\}.
\Ee
\elem
\begin{proof}
Let us put $$A:=\{z\in \mathcal H: Mf(z)>\lambda\}$$ and $$B:=\{z\in \mathcal H: M_df(z)>\frac{\lambda}{68}\}.$$ Recall that there is a family $\{Q_{I_j}\}_{j\in \mathbb N_0}$ of maximal (with respect to the inclusion) disjoint dyadic Carleson boxes (i.e $I_j\in \mathcal D$) such that
$$\frac{4\lambda}{68}\ge \frac{1}{|Q_{I_j}|}\int_{Q_{I_j}}|f|dV>\frac{\lambda}{68}$$
so that $B=\cup_{j\in \mathbb N_0}Q_{I_j}$.
\vskip .1cm
Let $z\in A$ and suppose that $z\notin B$. We know that there is an interval $I$ (not necessarily dyadic) such that $z\in Q_I$ and
\Be\label{eq:hypt} \frac{1}{|Q_I|}\int_{Q_I}|f|dV>\lambda.\Ee
Recall with Lemma \ref{lem:dyacovering} that $I$ can be covered by at most two adjacent dyadic intervals $J_1$ and $J_2$ (in this order) such that $|I|<|J_1|=|J_2|\le 2|I|$ so that $Q_I\subset Q_{J_1}\cup Q_{J_2}$. Of course, $z$ belongs only to one (and only one) of the associated boxes $Q_{J_1}$ and $Q_{J_2}$. Let us suppose that $z\in Q_{J_1}$. Then necessarily, $Q_{J_1}$ is not contained in $B$ since if so, $z$ would belong to $B$ and this would contradict our hypothesis on $z$. Thus $Q_{J_1}\cap B=\emptyset$ or $Q_{J_1}\supset Q_{I_j}$ for some $j$ and in both cases, because of the maximality of the $I_j$s, we deduce that
$$\frac{1}{|Q_{J_1}|}\int_{Q_{J_1}}|f|dV\le \frac{\lambda}{68}.$$
For the other interval $J_2$, we have the following possibilities
$$
\left\{ \begin{matrix} J_2=I_j\,\,\,\textrm{for some}\,\,\,j\\

               J_2\subset I_j \,\,\,\textrm{for some}\,\,\,j\\
               J_2\supset I_j\,\,\,\textrm{for some}\,\,\,j\\
               J_2\cap B=\emptyset.
                              \end{matrix}\right.
$$
If $J_2\supset I_j$ for some $j$ or $J_2\cap B=\emptyset$, then because of the maximality of the $I_j$s,

$$\frac{1}{|Q_{J_2}|}\int_{Q_{J_2}}|f|dV\le \frac{\lambda}{68}.$$

If $J_2=I_j$ for some $j$, then of course,
$$\frac{1}{|Q_{J_2}|}\int_{Q_{J_2}}|f|dV\le \frac{4\lambda}{68}.$$
It remains to consider the case where $J_2\subset I_j$ for some $j$. If $J_2\subset I_j$, then we can have
$$
\left\{ \begin{matrix} J_2=I_j^-\\

               J_2\subset I_j^- \\
               J_2\subseteq I_j^+
                              \end{matrix}\right.
$$
where $I_j^-$ and $I_j^+$ denote the left and right halfs of $I_j$ respectively. If $J_2\subset I_j^-$ or $J_2\subseteq I_j^+$, then $J_1\cap I_j\neq \emptyset$, and this necessarily implies that $J_1\subset I_j$. Thus $z\in Q_{J_1}\subset Q_{I_j}\subset B$ which contradicts the hypothesis $z\notin B$. Thus the only possible case is $J_2=I_j^-$ which leads to the estimate
$$\frac{1}{|Q_{J_2}|}\int_{Q_{J_2}}|f|dV\le \frac{4}{|Q_{I_j}|}\int_{Q_{I_j}}|f|dV\le \frac{16\lambda}{68}.$$

Thus from all the above analysis, we obtain
\Beas
\frac{1}{|Q_{I}|}\int_{Q_{I}}|f|dV &=& \frac{1}{|Q_{I}|}\left(\int_{Q_{I}\cap Q_{J_1}}|f|dV+\int_{Q_{I}\cap Q_{J_2}}|f|dV\right)\\ &\le&
\frac{|Q_{J_1}|}{|Q_{I}|}\left(\frac{1}{|Q_{J_1}|}\int_{Q_{J_1}}|f|dV+\frac{1}{|Q_{J_2}|}\int_{Q_{J_2}}|f|dV\right)\\ &\le& 4\left(\frac{\lambda}{68}+\frac{16\lambda}{68}\right)=\lambda
\Eeas
which clearly contradicts (\ref{eq:hypt}). The proof is complete.
\end{proof}
We can now prove Theorem \ref{thm:main3}.
\begin{proof}[Proof of Theorem \ref{thm:main3}.]
Let us note that by Lemma \ref{lem:equivdefBpq}, $(\textrm{b})\Leftrightarrow (\textrm{c})$. Let us prove that $(\textrm{a})\Leftrightarrow (\textrm{b})$.

Let $f$ be a locally integrable function and $I$ an interval. Fix $\lambda$ such that $0<\lambda<\frac{1}{|Q_I|}\int_{Q_I}|f|dV$. Then $$Q_I\subset \{z\in \mathcal H:M(\chi_{Q_I}f)>\lambda)\}.$$
It follows from the latter and (\ref{eq:main31}) that $$\mu(Q_I)\le \frac{C}{\lambda^q}\left(\int_{Q_I}|f(z)|^p\omega(z)dV(z)\right)^{q/p}.$$
As this happens for all $\lambda>0$, it follows in particular that
$$\mu(Q_I)\left(\frac{1}{Q_I}\int_{Q_I}|f|dV(z)\right)^q\le C\left(\int_{Q_I}|f(z)|^p\omega(z)dV(z)\right)^{q/p}.$$
\vskip .1cm
Next suppose that (\ref{eq:main33}) holds. We observe with Lemma \ref{lem:levelsetsembed} that to obtain (\ref{eq:main31}), we only have to prove the following
\Be\label{eq:main34}
\mu\left(\{z\in \mathcal D: M_df(z)>\frac{\lambda}{68}\}\right)\le \frac{C}{\lambda^q}\|f\|_{p,\omega}^q.
\Ee
We recall that $$\{z\in \mathcal H: M_df(z)>\frac{\lambda}{68}\}=\cup_{j\in \mathbb N_0}Q_{I_j}$$ where the $I_j$s are maximal dyadic intervals with respect to the inclusion and such that $$\frac{1}{|Q_{I_j}|}\int_{Q_{I_j}}|f|dV>\frac{\lambda}{68}.$$
Our hypothesis provides in particular that
$$\mu(Q_{I_j})\lesssim \left(\frac{|Q_{I_j}|}{\int_{Q_{I_j}}|f|dV}\right)^q\left(\int_{Q_{I_j}}|f|^p\omega dV\right)^{q/p}.$$
Thus
\Beas
\mu\left(\{z\in \mathbb D: M_df(z)>\frac{\lambda}{68}\}\right) &=& \sum_{j}\mu(Q_{I_j})\\ &\le& \sum_{j}\left(\frac{|Q_{I_j}|}{\int_{Q_{I_j}}|f|dV}\right)^q\left(\int_{Q_{I_j}}|f|^p\omega dV\right)^{q/p}\\ &\le& \left(\frac{68}{\lambda}\right)^q\sum_{j}\left(\int_{Q_{I_j}}|f|^p\omega dV\right)^{q/p}\\ &\le& \left(\frac{68}{\lambda}\right)^q\left(\sum_{j}\int_{Q_{I_j}}|f|^p\omega dV\right)^{q/p}\\ &\le& \left(\frac{68}{\lambda}\right)^q\left(\int_{\mathcal H}|f(z)|^p\omega(z) dV(z)\right)^{q/p}\\ &=& \left(\frac{68}{\lambda}\right)^q\|f\|_{p,\omega}^q.
\Eeas
The proof is complete.
\end{proof}
Taking $d\mu(z)=\sigma(z)dV(z)$, we obtain the following corollary.
\bcor\label{cor:main4}
Let $1\le p, q<\infty$, and $\omega,\sigma$ two weights on $\mathcal H$. Then the following assertions are equivalent.
\begin{itemize}
\item[(a)] There is a constant $C_1>0$ such that for any $f\in L_\omega^p(\mathcal H)$, and any $\lambda>0$,
\Be\label{eq:main41}
\sigma\left(\{z\in \mathcal H: Mf(z)>\lambda\}\right)\le \frac{C_1}{\lambda^q}\left(\int_{\mathcal H}|f(z)|^p\omega(z)dV(z)\right)^{q/p}
\Ee
\item[(b)] There is a constant $C_2>0$ such that for any interval $I\subset \mathbb R$,
\Be\label{eq:main42}
|Q_I|^{1/q-1/p}\left(\frac{1}{|Q_I|}\int_{Q_I}\omega^{1-p'}(z)dV(z)\right)^{1/p'}\left(\frac{1}{|Q_I|}\int_{Q_I}\sigma(z)dV(z)\right)^{1/q}\le C_1
\Ee
where $\left(\frac{1}{|Q_I|}\int_{Q_I}\omega^{1-p'}(z)dV(z)\right)^{1/p'}$ is understood as $\left(\inf_{Q_I}\omega\right)^{-1}$ when $p=1$.
\end{itemize}
\ecor

\bibliographystyle{plain}

\end{document}